\documentclass[10pt]{article}
\usepackage{amsmath}
\usepackage[colorlinks=true, pdfstartview=FitV, linkcolor=blue, citecolor=blue,urlcolor=blue]{hyperref}
\def\CSA{\mathfrak h}
\def\lie{\mathfrak g}

\usepackage{amsbsy}
\usepackage{amssymb}
\def\D{\mathbb D}

\usepackage{verbatim}
\def\M{\mathcal M}
\makeatletter
\@addtoreset{equation}{section}
\makeatother

\newtheorem{theorem}{Theorem}[section]

\newtheorem{exercise}{Exercise}[section]

\newtheorem{lemma}{Lemma}[section]
\newtheorem{remark}{Remark}[section]

\newtheorem{proposition}{Proposition}[section] 
\newtheorem{corollary}{Corollary}[section] 
\newtheorem{definition}{Definition}[section]
\newtheorem{conjecture}{Conjecture}
\def\le{\left}
\def\ri{\right}

\def\W{\mathfrak W}

\def\br{\begin{remark}}
\def\er{\end{remark}}
\def\bt{\begin{theorem}}
\def\et{\end{theorem}}
\def\k{\mathfrak k}
\def\bc{\begin{corollary}}
\def\ec{\end{corollary}}
\def\bx{\begin{examp}\small}
\def\ex{\end{examp}}
\def\bxr{\begin{exercise}\small}
\def\exr{\end{exercise}}
\def\bl{\begin{lemma}}
\def\el{\end{lemma}}
\def\bd{\begin{definition}}
\def\ed{\end{definition}}
\def\bp{\begin{proposition}}
\def\ep{\end{proposition}}
\def\be{\begin{equation}}
\def\ee{\end{equation}}
\def\ov {\overline}

\def\&{\hspace{-15pt}&}
\def\bea{\begin{eqnarray}}
\def\eea{\end{eqnarray}}
\def\beas{\begin{eqnarray*}}
\def\eeas{\end{eqnarray*}}
\def\B#1#2{\le\langle #1, #2 \ri \rangle}

\def \pa{\partial}
\def\C{{\mathbb C}}
\def\L{\mathcal L}
\def\R{{\mathbb R}}

\def\wh{\widehat}

\def\n{\mathfrak n}
\def\a{\alpha}
\def\d{\,\mathrm d}

\def\1{{\bf 1}}
\def\wt{\widetilde}
\newcommand{\tops}[2]{\texorpdfstring{#1}{#2}}
\def\omit#1{{}}

\date{}
\begin{document}
\baselineskip 15pt plus 1pt minus 1pt

\vspace{0.2cm}
\begin{center}
\begin{Large}
\fontfamily{cmss}
\fontsize{17pt}{27pt}
\selectfont
\textbf{Harish-Chandra integrals as  nilpotent integrals}
\end{Large}\\
\bigskip
\begin{large} {M.
Bertola}$^{\ddagger,\sharp}$\footnote{Work supported in part by the Natural
    Sciences and Engineering Research Council of Canada
(NSERC).}\footnote{bertola@crm.umontreal.ca}, A. Prats Ferrer$^{\sharp}$\footnote{pratsferrer@crm.umontreal.ca}.
\end{large}
\\
\bigskip
\begin{small}
$^{\ddagger}$ {\em Department of Mathematics and
Statistics, Concordia University\\ 1455 de Maisonneuve W., Montr\'eal, Qu\'ebec,
Canada H3G 1M8} \\ $^{\sharp}$ {\em Centre de recherches math\'ematiques, Universit\'e\ de Montr\'eal\\ 2920 Chemin de la tour, Montr\'eal, Qu\'ebec, Canada H3T 1J4.
}
\end{small}

\bigskip
{\bf Abstract}
\end{center}
Recently the correlation functions of the so--called Itzykson-Zuber/Harish-Chandra integrals were computed (by one of the authors and collaborators) for all classical groups using an integration formula that relates integrals over compact groups with respect to the Haar measure and Gaussian integrals over a maximal nilpotent Lie subalgebra of their complexification. Since the integration formula  {\em a posteriori} had the same form for the classical series, a conjecture was formulated that such a formula should hold for arbitrary semisimple Lie groups. We prove this conjecture using an abstract Lie--theoretic approach.

\vspace{0.7cm}

\section{Introduction and setting}
In random matrix theory \cite{DiFGiZJ, DaulKazKos} a particularly important role is played by the so--called Itzykson--Zuber/Harish-Chandra {measure. Its integral is} 
\cite{Harish-Chandra, IZintegral}
\be
\int_{U(N)} \d U {\rm e}^{tr (XUYU^\dagger)} = C_N \frac {\det({\rm e}^{x_i y_j} )_{i,j}}{\Delta(X)\Delta(Y)}
\ee
where $\d U$ is the Haar measure on $U(N)$, $X, Y$ are diagonal matrices and $\Delta(X)= \prod_{i<j} (x_i-x_j)$, with $X=diag(x_1,\dots, x_N)$ and $C_N$ is some proportionality constant. 

In a recent paper \cite{prats} {integrals of invariant functions on such type measures}
were computed for all the classical series ($A_n, B_n, C_n,D_n$) corresponding to unitary, orthogonal and symplectic ensembles, generalizing an earlier result \cite{prats0}.
An outstanding conjecture was formulated in \cite{prats} relating certain spherical integrals
\begin{conjecture}[Conjecture 1.1 in \cite{prats}]
\label{one}
Let $\lie$ a semisimple Lie algebra over $\C$, $\CSA$ a Cartan subalgebra, $K$ the maximal compact group in $\exp(\lie)$, $\B{}{}$ the Killling form and $\W$ the Weyl group. Let $F(X,Y)$ be  a polynomial on $\lie\times \lie$ invariant under diagonal adjoint action $F(Ad_g (X),Ad_g(Y)) = F(X,Y)$. Then the following identity holds $\forall H, J \in \CSA$:
\be
 \int_K \d k F(H, Ad_k J ){\rm e}^{-\B{H}{ Ad_k J}}  = \frac{C}{|\W|} \sum_{w \in \W} \epsilon_w \frac{{\rm e}^{-\B{H }{ J_w}}}{\prod_{\a>0} \a(H)\a(J)} \int_{\mathfrak n_+} \d N F(H + N , J_w  + N^\dagger) {\rm e}^{-\B{N}{N^\dagger}}\ ,\
\ee
where  $J_{w}$ stands for the action of the Weyl group on $\CSA$, $\epsilon_w$ is the usual sign homomorphism and $C$ is a  suitable constant depending only on the Lie algebra under consideration. (All the symbols will be defined more in detail later)
\end{conjecture}
Such conjecture was verified a posteriori for all the classical series but the authors of \cite{prats} failed to provide a general proof that would apply also to exceptional Lie algebras.  This identity was the main initial step towards an effective computation of all correlation functions for spherical integrals over the compact forms of the classical groups, $SU(N), SO(2n,\R), SO(2n+1,\R), Sp(n,\R)$. 
In this short note we provide a Lie--algebro--theoretical proof of Conjecture \ref{one} (Thm. \ref{mainthm}) that does not rely on any specificity of the Lie algebra as long as it is semisimple and provides a precise value for the proportionality constant $C$.

We will need to prove a slight generalization of Weyl--integration formula (which may well be known in the literature but we could not find in any of the standard references). The proof of Conjecture \ref{one} is contained in Thm. \ref{mainthm}.
%

We will  liberally use known facts about (semi)-simple Lie algebras, all of which can be found in standard reference books like \cite{Samelson}.
Let $\lie$ be a complex semisimple Lie algebra over $\C$, $G$ the corresponding simply connected group, $G = \exp(\lie)$.
Let $\CSA = \CSA_\C\subset \lie$ be a Cartan subalgebra (over $\C$ unless otherwise specified) and $\mathfrak R\subset \CSA_\C^\vee$ be the set of roots. Let an ordering of the roots be chosen: it fixes the set of positive roots $\mathfrak R_+$. The set of simple positive roots with respect to this ordering will be denoted by  $\Phi$.
 We will use the Chevalley basis $\{E_\alpha,H_\alpha\}_{\a\in \mathfrak R}$ (called {\em root vectors} and {\em coroots} respectively) of $\lie$  where $E_\alpha$ spans $\lie_\alpha$ and the set $H_\alpha$ for $\alpha\in\Phi$ spans $\CSA$. Such a basis has the properties $
[E_\alpha,E_{-\alpha}]= H_\alpha\ ,\ [H_\alpha, E_{\pm\alpha}] =\pm 2E_{\pm\alpha}\ ,\ \a\in \mathfrak R$.
Here and in the following we will use the notations 
\bea
\CSA_\R:= \sum_{\alpha \in  \Phi} \R\{H_\alpha\}\ ,
\qquad \CSA:= \CSA_\C:= \sum_{\alpha\in \Phi} \C\{H_\alpha\}\ ,\ 
\mathfrak n_+ =  \sum_{\alpha>0} \C \{E_{\alpha}\}\ ,\qquad
\mathfrak b_+  = \CSA_\C + \mathfrak n_+\ .
\eea
The compact form $\k\subset \lie$ will be chosen as the span of 
\be
\k = i\CSA_\R + \sum_{\a>0} \R\{X_\a, Y_\a\}\ ,\ \ X_\a:= (E_\a-E_{-\a})\ ,\qquad Y_\a:=i(E_\a+E_{-\a})
\ee
On each of the above Lie algebras we will use the Lebesgue measure such that the unit cube in the coordinates given by the specified basis has unit volume.


We finally recall that any semisimple Lie algebra  admits a decomposition $
\lie = \mathfrak k + \mathfrak a + \mathfrak n_+
$
where $\mathfrak k $ is a compact Lie algebra, $\mathfrak a\subset \CSA_\R$ is Abelian. This decomposition subtends the {\em Iwasawa decomposition}
of the Lie group $ G = K\,A\,N$.

\section{Schur decompositions}\label{Sec:dec}

The first fact we need is a generalization of the Schur decomposition to arbitrary Lie algebras: Schur decomposition is a widely known decomposition of matrices and states that any complex square matrix $M$ can be written in the form $U T U^\dagger$ with $U\in U(N)$ and $T$ an upper semi-triangular matrix.

\bt
\label{schurregular}
Any {\em ad-regular} element $M \in \lie$ is $K$--conjugate to an element in ${\mathfrak b}_+ = \CSA_\C + \mathfrak n_+$.
\et
{\bf Proof.}
It is known \cite{Samelson} that any ad-regular element in $\lie$ is conjugate to an element in $\CSA$, namely $ M = Ad_g H \in \CSA_\C$. 
There are (generically) $|\W|$ such ways of representing $M$.
Using the Iwasawa decomposition $g = k a n$ we have immediately
\be
M = Ad_k (H + (Ad_{an} H-H)) = Ad_k B
\ee
with $B = Ad_{an} H\in \mathfrak b_+$ and $Ad_{an} H-H\in \n_+$. {\bf Q.E.D.}\par\vskip 5pt
Since we will be concerned with integration formul\ae, the above theorem suffices since the set of ad-regular elements is a Zariski open set, dense in $\lie$: in particular nonregular elements are a set of Lebesgue measure zero. However, the following more general theorem can also be proved (but we will not prove it here in the interest of conciseness and also because it is completely irrelevant to our main purpose).

\bt
\label{schur}
Any element $M\in \lie$ is conjugated by an element of the maximal compact subgroup $K$ to an element $H+N$ with $H\in \CSA$, $N\in \mathfrak n_+$
\et

\section{Complex Weyl integration formula}
\label{CWIF}
The goal of this section is to write an integral formula for functions on $\lie$ in terms of integrals on $\mathfrak b$ and $K$.

Define 
\be
\M:= (K\times \mathfrak b)/T
\ee
where the action of the Cartan torus is $
t\cdot (k,V):= (kt^{-1}, tVt^{-1})$, $V\in \mathfrak b\ k\in K$ .
We will write $V = H + N$ with $iH\in \CSA_\C$ and $N\in \n_+$.  The action of $T$ is then 
$
t\cdot (k,H, N) = (kt^{-1}, H, t N t^{-1})
$
so that we can also think of $\M$ as $\M = \CSA \times (K\times \n_+)/T$.
The tangent space to $\M$ at $[(k,V)]$ is identified with $\mathfrak k/i\CSA_\R + \mathfrak b$ by $[(k_s,V_s)]:= [( k{\rm e}^{sX},V+sW)]$.

Consider the map 
\be
\begin{array}{rcl}
\pi:\M & \longrightarrow & \lie\\[1pt]
[k,V] & \mapsto k V k^{-1}
\end{array}
\ee
The topological degree of $\pi$ is the cardinality of the Weyl group $\W$: to see this it is sufficient to note that $Ad_K \CSA_\C \subset \lie$ has the advocated degree  and then use a continuity argument.

The differential of the map $\pi:\M \to \lie$ at a point $[(k,V)]$ is then computed as 
\be
\frac \d{\d s} Ad_{k{\rm e}^{sX}} (V+sW)\big|_{s=0} = k\le([X,V] + W\ri)k^{-1}\ ,\ \ X\in \k/i\CSA_\R,\ W\in \mathfrak b\ .
\ee
In order to write a matrix representation of the above map  we write it in the natural basis of $T_{M}\lie = \n_- + \mathfrak b_+$
\bea
\begin{array}{rcl}
d\pi: T_{[k,V]}\M\sim \mathfrak k/i\CSA_\R \oplus \mathfrak b_+ &\longrightarrow& \n_- + \mathfrak b_+\\
(X,W)& \mapsto& \d\pi(X,W) = Ad_k([X,V] +W)
\end{array}
\eea
We compute the determinant of the above map without the $Ad_k$ term (which does not change its value) and we think of $\n_-$ as a vector space over $\R$ with a real basis provided by $ \n_-:= \sum_{\a<0} \R\{E_\a\} + i \sum_{\a<0} \R\{E_\a\}$
\bea
ad_{H+N}(X_\a) = \Re(\a(H)) E_{-\a} + i \Im(\a(H)) E_{-\a}+  \sum_{-\beta>-\a} \C\{E_{-\beta}\}\  {\rm mod}\  \mathfrak b_+\\
ad_{H+N}(Y_\a) = -\Im(\a(H)) E_{-\a} +i \Re(\a(H))E_{-\a}+  \sum_{-\beta>-\a} \C\{E_{-\beta}\} \ {\rm mod}\ \mathfrak b_+\\
V:= H+N\nonumber
\eea
It appears that the matrix has a block-uppertriangular shape and  these ``upper triangular'' parts  do not contribute to the determinant. The latter becomes then the product of the determinants of the above $2\times 2$ blocks which are simply $|\a(H)|^2$.

The Jacobian of $\d\pi$ at the point $[k,V]\in \M$ ($V=H+N$) is thus 
\be
J(k,V) = \prod_{\alpha\in\mathfrak R_+} \le|\alpha(H)\ri|^2=: |\Delta(H)|^2.
\ee
The notation $\Delta(H):= \prod_{\alpha\in\mathfrak R_+} \alpha(H)$ is used in analogy with the case of $\lie  = sl(n,\C)$ where it reduces to the Vandermonde determinant. A well known property is that 
\be
\Delta(H_w) = (-)^w \Delta(H)
\ee
where $w\in \W$ and $H_w$ stands for the action of the Weyl group on $\CSA_\C$ and the notation $(-)^w$ means the parity of the Weyl-transformation (i.e. the parity of the number of elementary reflections  along walls of Weyl chambers in which $w$ can be decomposed).
Collecting these pieces of information we have proved the following
\bt[Complex-Weyl integration formula]
\label{complexweyl}
Let $F:\lie\to \C$ be a smooth integrable function invariant under the adjoint action of $K$. Then 
\bea
\int_\lie \d M F(M)  =c_\k \int_{\CSA_\C \times \mathfrak n_+} \hspace{-10pt}\d H \d N F(H+N) 
\le|\Delta(H)\ri| ^2
\label{Jaco.1}\ ,\qquad c_\k:= \frac {\mu(K)/\mu(T)}{|\W|}
\nonumber
\eea
where $\d M,\ \d H,\ \d N$ are the Lebesgue measures on  $\lie,\ \CSA_\C,\ \mathfrak n_+$ respectively defined above,  $\mu(K)$ and $\mu(T)$ are the induced measures on the compact group $K$ and the maximal torus $T$, $\W$ is the Weyl group and $|\W|$ is its cardinality.
\et

There is one more piece of information that we can extract from the above and is contained in the following 
\bc
\label{Winv}
For any $Ad_K$--invariant smooth integrable function $F:\lie\to\C$ the function 
\be
\wh F(H):= \int_{\n_+} \d N F(H+N): \CSA_\C\to \C
\ee
is Weyl--invariant.
\ec
{\bf Proof.}
By the generalized Schur decomposition (Thm. \ref{schurregular}) a regular element $M$ can be represented modulo the $Ad_K$ action as $H+N\in \CSA + \n_+$ or $H_w + \wt N$ where $H_w$ is in the same $\W$--orbit through $H$ and $\wt N\in \n_+$ is some other element in the {\em same} nilpotent subalgebra $\n_+$. In general the dependence of $\wt N$ on $N, H$ is  a complicated expression. Consider a small ball  $M\in \mathcal U\subset \lie$ consisting of regular elements. This ball can be mapped {\em diffeomorphically} to some neighborhood $\mathcal H \times \mathcal L$ in $\CSA_\C \times (K\times \n_+)/T$ with $\mathcal H$ lying in a suitable Weyl chamber and containing $H$. Choosing another Weyl chamber $\mathcal H_w$ we have a distinct diffeomorphism between $\mathcal U$ and $\mathcal H_w\times \wt {\mathcal L}$. Since the Jacobian computed above is always $|\Delta(H)|^2 = |\Delta(H_w)^2|$ and {\em independent} of $N$, we conclude that the transformation $N\mapsto \wt N$ preserves the Lebesgue measure of $\n_+$, namely $\d \wt N/\d N=1$. Since regular elements are open and dense (and the complement has zero measure) we can then write 
\be
\wh F(H_w) = \int_{\n_+} \d \wt N F(H_w+\wt N) = \int _{\n_+}\d \wt N F(Ad_k (H_w+\wt N)) = \int _{\n_+} F(H+N)\d N = \wh F(H)\ .
\ee
{\bf Q.E.D.}\par \vskip 5pt

The integration formula of Thm. \ref{complexweyl} should be considered a mild generalization of the standard Weyl integration formula which we state here for functions on the Lie algebra of a compact Lie group $K$.

\bt[Weyl integration formula]
\label{weyl}
Let $F:\k\to \C$ be a smooth function integrable with respect to the Lebesgue measure. Then 
\be
\int_\k \d X F(X) = c_\k  \int_{i\CSA_\R} \d H |\Delta(H)|^2 \int_{K} \d k F(Ad_k(H))
\label{Jaco.2}
\ee
In case of an $Ad_K$--invariant function the above reduces to 
\be
\int_\k \d X F(X) = c_\k \int_{i\CSA_\R} \d H |\Delta(H)|^2  F(H)
\ee
where we  need to put the absolute--value sign because $\Delta(H)^2$ may be negative on $i\CSA_\R$  if $\dim n_+$ is odd ($c_\k$ has the same meaning and value as in Thm. \ref{complexweyl}).
\et
\br
The explicit value of $\mu(K)/\mu(T)$ was computed by Macdonald  \cite{Macdonald} for a slightly different choice of normalization for the Lebesgue measure. The actual value of this constant is irrelevant for our purpose and does not directly enter the computation of the proportionality constant in the conjecture.
\er
\section{Gaussian integrals}
Let $V_{{\R}}$ be a real vector space and $\B{}{}:V_{{\R}}\times V_{{\R}}\to \R$ a positive definite bilinear pairing (an inner product).
Let $V_\C:= V_{{\R}}\otimes \C$ be its complexification.  Denote by $\d x$ the Lebesgue measure on $V_{{\R}}$ and $\d z$ the Lebesgue measure on $V_\C$ (the normalizations of which are irrelevant at this point). 
The inner product $\B{}{}$ extends to an inner-product (linearly over $\C$) on $V_\C$. Moreover the real form $V_\R\subset V_\C$ defines also a natural conjugation $z\to \ov z$ which fixes $V_\R$. 

\bl 
\label{LemmaWick}
With the notations and definitions above, for any polynomial function $F$ on $V_{{\C}} \times V_{{\C}}$ define
\bea
&& <F>_\R:= \frac 1 {Z_\R} \int_{V_{{\R}}}\int _{V_{{\R}}}\d x \d y {\rm e}^{-a\B xx - c\B yy -2b \B xy} F(x,y)\cr
&&  <F>_\C:= \frac 1{Z_\C} \int_{V_\C} \d z  {\rm e}^{-a\B z z - c\B {\ov z}{\ov z} -2b \B z{\ov z}} F(z,\ov z)
\eea
where $Z_\R$ and $Z_\C$ are determined\footnote{It is an easy exercise that we leave to the interested reader to verify that $\mathcal Z_\R =  (2\pi)^{n}  \delta^{  \frac n2}\ ,\
\mathcal Z_\C  = (2\pi)^{n} (-\delta)^{\frac n2}\ ,\ \ \ n=\dim_\R V_\R =\dim_\C V_\C ,\ \  \delta = ac-b^2.$ These precise expressions are nevertheless irrelevant for our purposes.
}
  by the requirement that $<1>=1$.
While the convergence of the integrals in the two cases imposes different conditions on the numbers $a,b,c$, nevertheless 
{\bf (i)} both $ <F>_{\R,\C}$ are polynomials in $a/\delta, b/\delta, c/\delta$, $\delta:= {ac -b^2}$ and {\bf (ii)} as polynomials they coincide.
\el
{\bf Proof}.
The key is in showing that the generating functions for the moments of the two integrals in either cases are identical, namely that for $A,B\in V_\C$
\bea
 G_\R(A,B):= \le<{\rm e}^{\B xA + \B y B}\ri>_{\R} = \exp \le[\frac 14 \le( \frac c \delta \B AA + \frac a \delta \B BB -2 \frac b\delta \B AB\ri)\ri] = G_\C(A, B):= \le<{\rm e}^{\B z A + \B {\ov z} B }\ri>_{\C} 
 \label{52}
\eea
for then 
$<F>_\R = F(\pa_A, \pa_B) G_\R(A,B)\big\vert_{A=0=B} = 
<F>_\C $, which proves both points of the lemma at the same time.
In order to show (\ref{52})  we use an orthonormal coordinate basis for $\B{}{}$ so that --writing $A= (\a_1,\dots, \a_n)$ and $B = (\beta_1,\dots, \beta_n)$ in this basis-- the integral $G_\R$ factorizes as $G_\R = \prod G_1(\a_j,\beta_j)$ with  
\be
G_1(\a,\beta):= \frac 1{Z_{1,\R}}\int_{\R} \d x \int _\R\d y\, {\exp} \le[-(x,y) M \begin{pmatrix}x\cr y\end{pmatrix} + x\a + y \beta\ri]\ ,\ \ M:=\begin{pmatrix} a&b\cr b & c\end{pmatrix} 
\ee
Define $(x',y') = (x,y) - \frac  12 (\a,\beta)M^{-1}$, where now the contours of integration may be some lines parallel to the real axis in the complex $x'$ and $y'$ planes. However the ensuing integrals can be deformed (by Cauchy theorem) back to the real axis and the integral yields $G_1(\a,\beta) = g(\a,\beta):= \exp \frac 14 \le[\frac c\delta \a^2 + \frac a \delta \beta^2  - 2\frac b \delta \a\beta\ri]$.

For the second case we have $G_\C =\prod \wt G_1(\a_j,\beta_j)$ with $\wt G_1(\a,\beta)$ given below:  we need to express the integration in the real/imaginary part of $z = x+iy$
\bea
\wt G_1(\a,\beta):=\frac 1{Z_{1,\C}}\int_{\C} \d^2z \, {\exp} \le[-(z,\ov z)\begin{pmatrix} a&b\cr b & c\end{pmatrix} \begin{pmatrix}z\cr \ov z\end{pmatrix} + z\a + \ov z \beta\ri] =\\
=\int_{\R} \int_\R\d x\d y\,  {\exp} \le[-(x,y)\begin{pmatrix} a+c+2b &i(a-c)\cr i(a-c)  & 2b -a -c\end{pmatrix} \begin{pmatrix} x \cr y\end{pmatrix} + x(\a + \beta) + i y(\a - \beta)\ri]  
\eea
In this case we perform the shift $(x',y')=(x,y) -\frac 1 2 (\a+\beta,i(\a-\beta)) \begin{pmatrix} a+c+2b &i(a-c)\cr i(a-c)  & 2b -a -c\end{pmatrix}^{-1}$, followed by deforming back the integration contours on the real $x'$ and $y'$ axes. Straightforward linear algebra gives the same  result $g(\a,\beta)$ as above.
{\bf Q.E.D.}\par \vskip 5pt


\section{Proof of Conjecture \tops{\ref{one}}{1}}\label{Sec:TMT}

Let $\vartheta: \lie\to\lie$ be the {\bf Cartan involution} (antilinear)
\be
\vartheta(cE_\a)= -\overline cE_{-\a}\ ,\ \ \vartheta(c H_\alpha) = -\overline c H_\alpha.
\ee

\bd
For a (semi)simple Lie algebra $\lie$ over $\C$, given the Cartan involution $\vartheta$ defined above, we will denote by $M^\dagger = -\vartheta(M)$, and by $M^\vartheta = \vartheta(M)$.
\ed
\br
The notation $M^\dagger$ has been defined to coincide with the usual hermitian conjugate in the standard fundamental representation of $sl(n,\C)$.
\er

The Cartan involution $\vartheta$ (or $\dagger$)  defines a real form of $\lie$ which is precisely $\k$, the {\em compact real form} as the $1$--eigenspace of $\vartheta$. Two properties are immediate ($\B{}{}$ is the Killing form):
\begin{itemize}
\item $\B{M}{M^\vartheta}\leq 0$ is a {\bf negative} definite {\bf sesquilinear} quadratic form for  $M\in \lie$;
\item $\B{X}{X^\dagger}=-\B{X}{X^\vartheta}=-\B{X}{X}\geq 0$ is  a {\em positive} definite quadratic form for $X\in \k$ (as a {\bf real} vector space).
\end{itemize}


Consider the following quadratic form on $\lie\times \lie$ 
\be
Q_A(X,Y):= a \B{X}{X} + 2b \B{X}{Y} + c \B{Y}{Y}\ ,\qquad A:= \begin{pmatrix} a&b\cr b&c \end{pmatrix}
\ee
We leave to the reader to verify the following easy
\bl
There exist two open domains $\D_{\k}$ and $\D_{\i\lie}$ for the parameters $a,b,c$ such that
\begin{itemize}
\item if $(a,b,c)\in \D_{\k}$ then  $\Re\, Q_A(X,Y)$ is positive definite on $\mathfrak k\times \mathfrak k$ ;
\item if $(a,b,c)\in \D_{\i\lie}$ then $\Re\,Q_A (M,M^\vartheta)$ is  positive definite on $\lie$.
\end{itemize}
\el
The specific form of these domains is largely irrelevant for our considerations and in the interest of conciseness we will not specify them further.
\bd
Define 
\bea
&& Z_\k:= \int_\k\int_\k dX dY {\rm e}^{-Q_A(X,Y)}\ , \ \ (a,b,c)\in \D_{\k} \ ,\qquad 
 Z_\lie:= \int_\lie dM {\rm e}^{-Q_A(M,M^\vartheta)}\ , \ \ (a,b,c)\in \D_{\i\lie}\ .
\eea
Then, for any {\em polynomial} function $F$ on $\lie\times \lie$ we define
\bea
<F>_\k:= \frac 1{Z_\k} \int_\k\int_\k dX dY F(X,Y) {\rm e}^{-Q_A(X,Y)}\ ,
\qquad
<F>_\lie:=\frac 1{Z_\lie} \int_\lie dM F(M,M^\vartheta) {\rm e}^{-Q_A(M,M^\vartheta )}
\eea
\ed

As an application of Lemma \ref{LemmaWick} with $V_{{\R}} =\k$, $V_\C = \k \otimes \C = \lie$ and $\vartheta$ as the involution leaving $\k$ invariant, we have 
\bp
\label{PropWick1}
For any polynomial function $F$ on $\lie\times \lie$ both $<F>_\k$ and $<F>_\lie$ are {\bf polynomials} in $a/\delta,b/\delta,c/\delta$ with $\delta = ac-b^2$. As polynomials they {\bf coincide}
\label{prop:GI}
\ep

The next theorem contains the proof of Conjecture \ref{one} with precise values of the proportionality constants.
\bt
\label{mainthm}
Let $F$ be an $Ad_K$ invariant (polynomial) function on $\lie\times \lie$,  where the action of $Ad_K$ is the diagonal action
\be
F(X,Y) = F( Ad_k(X), Ad_k(Y))\ ,\ \ \forall X,Y\in \lie
\ee
Then, for any $H,J\in \CSA_\C$
\bea
\int_K \d k F(H,Ad_k(J)) {\rm e}^{
\gamma\B{H}{Ad_k(J)}} 
= \frac {C_\lie}{|\mathfrak W|} 
\sum_{w\in\mathfrak W} \frac{{\rm e}^{
\gamma\B{H}{J_w}}} {\Delta(H)\Delta(J_w)} 
\frac{\int_{\mathfrak n_+} \d N F(H+N, J_w+N^\vartheta) {\rm e}^{
\gamma\B{N}{N^\vartheta}}}
{\int_{\mathfrak n_+} \d N {\rm e}^{
\gamma\B{N}{N^\vartheta}}}
\eea
The normalization constant $C_\lie$ is given by 
\be
C_\lie =
 |\W| \prod_{j=1}^{\dim \CSA} m_j! \prod_{\a>0}\frac{\B{\a}{\a}}{
2\gamma}\ .
\ee
where $m_j$ are the exponents of the Weyl group.
{Moreover  we have
\be
\int_{\mathfrak n_+} \d N {\rm e}^{
\gamma\B{N}{N^\vartheta}} =\prod_{\a>0} \frac {\pi \B \a\a} {
2\gamma}
\label{Zn}
\ee}
\et
\br
By simplifying the values of the constants (note that $\dim \n_+$ is the number of positive roots) we obtain 
\be
\int_K \d k F(H,Ad_k(J)) {\rm e}^{
\gamma\B{H}{Ad_k(J)}} 
= 
 \frac{\prod_{j=1}^{\dim \CSA} m_j!}{\pi^{\dim \n_+}}   
\sum_{w\in\mathfrak W} \frac{{\rm e}^{
\gamma\B{H}{J_w}}}{\Delta(H)\Delta(J_w)} 
{\int_{\mathfrak n_+} \d N F(H+N, J_w+N^\vartheta) {\rm e}^{
\gamma\B{N}{N^\vartheta}}}
\ee
\er
\br
{Note that the convergence of the Gaussian integral in the formula demands $\Re(\gamma)>0$, but the identity is one between analytic functions of $\gamma$.}
\er
{\bf Proof}.
We start from the proof of (\ref{Zn}): since $\B{E_\a}{E_\beta}=\frac{2}{\B \a\a}\delta_{\a,-\beta}$\footnote{ Indeed, $\B{[E_\a,E_{-\a}]}{ H } = \B{E_\a}{[E_{-\a},H]} = \a(H) \B{E_\a}{E_{-\a}}$. On the other hand $[E_\a,E_{-\a}] = H_\a$ and thus $\a(H) \B{E_\a}{E_{-\a}} = \B{H_\a}{H}$. Evaluating on $H=H_\a$ ($\a(H_\a)=2$, $\B{H_\a}{H_\a} = \frac 4{\B{\a}{\a}}$) we get the assertion.
}, writing
$N=\sum_{\a>0} n_\a E_\a$ the integral is recast into the form 
\be
\int_{\mathfrak n_+} \d N {\rm e}^{
\gamma\B{N}{N^\vartheta}} = \prod_{\a>0} \int_\C \d^2 n_\a {\rm e}^{\frac{-2\gamma}{\B \a\a} |n_\a|^2} = \prod_{\a>0}\frac {\pi \B \a\a} {
2\gamma}\ \ {(\Re(\gamma)>0)}\ .
\ee
The value of $C_\lie$ is computed by evaluating explicitly the integrals on both sides for $F\equiv 1$, which reduces the formula to the famous Harish--Chandra expression
\be
\int_K \d k {\rm e}^{
\gamma\B{H}{Ad_k(J)}} 
= \frac {C_\lie}{|\mathfrak W|} 
\sum_{w\in\mathfrak W} \frac{{\rm e}^{
\gamma\B{H}{J_w}}}{\Delta(H)\Delta(J_w)} 
\ee
 In this case the equality was established in  (Thm. 2, pag 104 \cite{Harish-Chandra})\footnote{In loc. cit. the exponent has a plus sign and no constant $\gamma$, which means that we have to map $H\mapsto 
{\gamma} H$ in Harish--Chandra's formula, thus yielding the {factor $(
\gamma)^{-\dim \n_+}$} due to the homogeneity of $\Delta$.}  where the value of the constant $C_\lie$ was given by $C_\lie =(
{\gamma})^{{-}\dim \n_+} \big(\Delta,\Delta\big)$.
 The bracket $\big(p(H),q(H)\big)$ was defined ibidem for any polynomials $p,q$ over $\CSA$ by writing them in a {\em orthonormal} basis $\B{\omega_j}{\omega_k}=\delta_{jk}$\footnote{We denote by the same symbol $\B{}{}$ the induced inner product on $\CSA^\vee$.}
\be
p(H):= \sum_{\vec n} a_{\vec n} \prod_{\ell=1}^{\dim(\CSA)}\omega_\ell^{n_\ell}(H)\  ,\ \  q(H):= \sum_{\vec n} b_{\vec n} \prod_{\ell=1}^{\dim(\CSA)}\omega_\ell^{n_\ell}(H)
\ee
\be
\big( p,q\big) := \sum_{\vec n} a_{\vec n}b_{\vec n} \prod_{\ell=1}^{\dim(\CSA)} m_\ell!
\ee

The number $\big(\Delta,\Delta\big)$ ($\Delta = \prod_{\a>0}\a$) has been computed in \cite{Macdonald}\footnote{ The formula is quoted as reported in an appendix of a paper of Harder cited ibidem, due to a private communication of Steinberg.} and is given by 
\be
\big(\Delta,\Delta\big) = {2^{\frac{\dim \CSA-\dim\lie}2}} |\W| \prod_{j=1}^{n} m_i! \prod_{\a>0} \B\a\a 
\ee
which proves the expression for the constants noticing that $\frac{\dim\lie-\dim\CSA}2= \dim \n_+$.

We now turn to the proof of the equality: consider first the integral over $\k\times \k$
\bea
\frac 1{\mathcal Z_\k}&\&\int_\k\!\int_\k \d X \d Y F(X,Y){\rm e}^{-a\B XX - c\B YY -2b \B XY}  = \cr
&\&= \frac {{c_\k^2}}{\mathcal Z_\k}  \int_{i\CSA_{_\R}} \int_{i\CSA_{_\R}} \d H \d J \Delta(H)^2 \Delta(J)^2 {\rm e}^{- a \B HH  - c \B JJ}
\underbrace{ \int_{K} \d k F( H, kJk^{-1})  {\rm e}^{- 2b \B H{kJk^{-1}}}}_{=:I(H,J)} \label{real}
\eea
where we have used Weyl integration formula (Thm. \ref{weyl}) twice. 
On the other hand for the integral over $\lie$,  using the complex Weyl integration formula (Thm. \ref{complexweyl})  we have 
\bea
 T(a,b,c)&\& :=\frac 1{\mathcal Z_\lie} \int_\lie  \d M F(M,M^\vartheta ){\rm e}^{-a\B MM - c\B {M^\vartheta}{M^\vartheta} -2b \B M {M^\vartheta}}  = \cr
&\& =\frac {{c_\k}}{\mathcal Z_\lie}  \int_{\CSA_{_\C}}  \d Z |\Delta(Z)|^2  {\rm e}^{-a \B ZZ -c \B {Z^\vartheta}{Z^\vartheta}} \int_{\mathfrak n_+} \d N F(Z+N, Z^\vartheta+ N^\vartheta )  {\rm e}^{-2b \B {Z+N}{Z^\vartheta + N^\vartheta}} =\cr
&\& =\frac {{c_\k}}{\mathcal Z_\lie}  \int_{\CSA_{_\C}}  \d Z |\Delta(Z)|^2  {\rm e}^{-a \B ZZ -c \B {Z^\vartheta}{Z^\vartheta} -2b \B Z{Z^\vartheta}} \int_{\mathfrak n_+} \d N F(Z+N, Z^\vartheta+ N^\vartheta )  {\rm e}^{-2b \B {N}{ N^\vartheta}} 
\eea
Here we have used that $f(M):= F(M,M^\vartheta){\rm e}^{-Q_A(M,M^\vartheta)}$ is $Ad_K$--invariant (but not $Ad_G$--invariant!) and then the simple fact that $\B{Z+N}{Z^\vartheta+N^\vartheta} = \B{Z}{Z^\vartheta} + \B{N}{N^\vartheta}$.

\noindent We now point out that $T(a,b,c)$ is a polynomial in $a/\delta,b/\delta,c/\delta$ by  Lemma \ref{LemmaWick}. Since $|\Delta(Z)|^2 = (-)^{\n_+}\Delta(Z)\Delta(Z^\vartheta)$ (where $Z^\vartheta=-Z^\dagger$ is again the natural conjugation w.r.t. $i\CSA_\R$), applying once more Lemma \ref{LemmaWick} with $V_{{\R}} = {i}\CSA_\R$ and $V_\C = \CSA_\C$ we obtain 
\bea
&\& T(a,b,c)  = \frac {(-)^{\n_+}{c_\k}}{\mathcal Z_\lie}  \int_{{i}\CSA_{_\R}}\!\!\!\d H\int_{{i}\CSA_{_\R}}\!\!\!  \d J\,   \Delta(H)\Delta (J)   {\rm e}^{-a \B HH -c \B JJ -2b \B HJ }  \int_{\mathfrak n_+} \d N F(H+N, J + N^\vartheta )  {\rm e}^{-2b \B {N}{ N^\vartheta}} = \cr
&\& = \frac {(-)^{\n_+}{c_\k}{\mathcal Z_{\n_+}}}{\mathcal Z_\lie}  \int_{{i}\CSA_{_\R}}\!\!\!\d H\int_{{i}\CSA_{_\R}}  \!\!\!\d J\,   \Delta(H)^2\Delta (J)^2   {\rm e}^{-a \B HH -c \B JJ}
\underbrace{\frac{{\rm e}^{ -2b \B HJ } }  {\Delta(H)\Delta(J)}  \frac{1}{{\mathcal Z_{\n_+}}}
\int_{\mathfrak n_+} \d N F(H+N, J + N^\vartheta )  {\rm e}^{-2b \B {N}{ N^\vartheta}}}_{=: G(H,J)}\cr 
&\&\label{imag}
\eea
where the last line is just a different way of rewriting the previous line with ${\mathcal Z_{\n_+}}=\int_{\n_+}\exp{-2b \B{N}{N^\vartheta}}\d N$.

Comparison of formul\ae\ (\ref{real}) and (\ref{imag}) suggests
the na\"ive observation  that 
$I(H,J) = \frac{(-)^{\n_+}{c_\k} {\mathcal Z}_\k{\mathcal Z}_{\n_+}}{{c_\k^2}\mathcal Z_\lie}  G(H,J)$ but this cannot possibly be the case since $I(H,J)$ is Weyl--invariant in {\bf both} variables while in general  $G(H,J)$ is not. 

What will be shown instead is that the {\em symmetrization} of $G(H,J)$ is proportional to $I(H, J)$, namely
\be
I(H,J) = \frac{(-)^{\n_+}{c_\k} {\mathcal Z}_\k {\mathcal Z}_{\n_+}}{{c_\k^2} \mathcal Z_\lie} {\frac{1}{|\W|}} \sum_{w\in \W}  G(H,J_w)\label{bingo}
\ee
which is precisely the assertion of our theorem.  {Note that $G(H_w,J_w)$: indeed from Corollary \ref{Winv} the integral\footnote{The prefactor of which already has the advocated invariance.}
\be
f(Z,Z^{{\vartheta}}) := \int_{\n_+}  \d N  F(Z + N, Z^\vartheta +N^\vartheta) {\rm e}^{-2b\B{N}{N^\vartheta}}
\ee
is a polynomial  (since $F$ is a polynomial in both variables) with invariance $f(Z_w, Z^\vartheta_w)= f(Z,Z^\vartheta)$. The polynomial $f(H,J)$ can be written as 
\be
f(H,J) =  {\rm e}^{H\pa_Z}{\rm e}^{ J \pa_{Z^\vartheta}} f(Z,Z^\vartheta)\big|_{Z=0=Z^\vartheta}
\ee
where $H\pa_Z, J\pa_{Z^\vartheta}$  stand for the vector--fields $H\pa_Z Z = H, H\pa_Z Z^\vartheta =0$ and similarly 
$J\pa_{Z^\vartheta} Z=0, \ J\pa _{Z^\vartheta} Z^\vartheta = J$.
}
Therefore it is sufficient to symmetrize  $G(H,J)$ with respect to --say-- $J$ in order to obtain a completely Weyl--invariant function. The symmetrization can be carried under the integral sign  without changing its value  since the measures $\d H, \d J$ and the exponential factors that precede $G$ in (\ref{imag}) are all $\W$--invariant. We have thus obtained
\be
 \frac {{c_\k}{\mathcal Z}_{\n_+}} {(-)^{\n_+}\mathcal Z_\lie} \int_{{i}\CSA_{_\R}\times i\CSA_\R}
 \hspace{-30pt}\d H\d J   \Delta(H)^2\Delta (J)^2   {\rm e}^{-a \B HH -c \B JJ} \sum_{w\in \W}\frac { G(H,J_w)}{|\W|} = \frac {{c_\k^2}}{\mathcal Z_\k}  \int_{{i}\CSA_{_\R}\times i\CSA_\R}\hspace{-30pt}\d H\d J
\Delta(H)^2 \Delta(J)^2 {\rm e}^{-a \B HH -c \B JJ}I(H,J)\nonumber
\ee

Of course this equality {\em per se} does not imply eq. \ref{bingo}. However we can use the following argument.
We replace the invariant polynomial $F(X,Y)$ by $h(X) g(Y) F(X,Y)$ with $h,g$ arbitrary $Ad_K$ invariant polynomials over $\k$

Note that $g(H+N)= g(H)$ (and so for $h$): indeed any $Ad_K$--invariant polynomial on $\k$ is automatically $Ad_G$--invariant  (on $\k\otimes \C = \lie$) and for a generic $H$, $H+N$ is $Ad_G$--conjugate to $H$ itself  since $ad_{H+N}$ is semisimple (in the adjoint representation).
  Thus, in eq. \ref{imag}, the two extra factors $h,g$ will be independent of $N$ and thus factorizable outside of the integral over the nilpotent algebra $\mathfrak n_+$.
  
  On the other hand, in eq. \ref{real} they clearly and immediately factor out of the $K$-integral thus yielding the identity
  \be
   \int_{{i}\CSA_{_\R}}\hspace{-10pt}\d H\int_{{i}\CSA_{_\R}} \hspace{-10pt}\d J\,  
   \Delta(H)^2\Delta (J)^2   {\rm e}^{-a \B HH -c \B JJ}h(H) g(J) 
 \underbrace{\left(\frac {{c_\k}{(-)^{\n_+}\mathcal Z}_{\n_+}} {\mathcal Z_\lie}{\frac{1}{|\W|}}\sum_{w\in \W} 		
 	G(H,J_w) -\frac {{c_\k^2}}{\mathcal Z_\k}I(H,J)\right)}_{=: R(J,H)} =0 \label{almostthere} 
  \ee
  valid for {\em arbitrary} Weyl--invariants polynomials $h,g$ on $\CSA_\R$.
 Note  that in $\mathcal H:= \L^2({i}\CSA_\R\times {i}\CSA_\R, \d H \d J {\rm e}^{-a\B HH - c \B JJ })$ the set of all polynomials  is dense and that  the bracket expression above belongs to this space (we can take $a,c\in \R_-$ for this computation). The projector onto the subspace of $\W$--invariant functions 
\be
\mathcal H_\W:=\{ f(H,J)\in \mathcal H:\ f(H_w,  J_{w'})  = f(H,J)\ ,\forall w,w'\in \W\}
\ee
is self-adjoint and hence the range is a closed subspace, to which $R(J,H)$ belongs. 
The space of Weyl invariant polynomials form a basis in this space and in particular are dense. Thus the vanishing of eq. \ref{almostthere} says that $R(J,H)$ is orthogonal to such a dense set, thus is identically vanishing. The last detail is that the identity so far has been proved only for $H,J\in i\CSA_\R$; however, being an identity between polynomials, it must hold for its complexification as well, namely on the whole $\CSA_\C$. The theorem is proved and so is Conjecture \ref{one}, { with $\gamma=-2b$}. {\bf Q.E.D.}\par \vskip 5pt

\bibliographystyle{unsrt}
\bibliography{Biblio}

\def\cydot{\leavevmode\raise.4ex\hbox{.}}
\begin{thebibliography}{1}

\bibitem{DiFGiZJ}
P.~Di~Francesco, P.~Ginsparg, and J.~Zinn-Justin.
\newblock {$2$}{D} gravity and random matrices.
\newblock {\em Phys. Rep.}, 254(1-2):133, 1995.

\bibitem{DaulKazKos}
J.-M. Daul, V.~A. Kazakov, and I.~K. Kostov.
\newblock Rational theories of {$2$}d gravity from the two-matrix model.
\newblock {\em Nuclear Phys. B}, 409(2):311--338, 1993.

\bibitem{Harish-Chandra}
Harish-Chandra.
\newblock Differential operators on a semisimple {L}ie algebra.
\newblock {\em Amer. J. Math.}, 79:87--120, 1957.

\bibitem{IZintegral}
C.~Itzykson and J.~B. Zuber.
\newblock The planar approximation. {II}.
\newblock {\em J. Math. Phys.}, 21(3):411--421, 1980.

\bibitem{prats}
A.~Prats~Ferrer, B.~Eynard, P.~Di~Francesco, and J.-B. Zuber.
\newblock Correlation functions of harish--chandra integrals over the
  orthogonal and the symplectic groups.
\newblock {\em J. Stat. Phys.}, 129(5-6):885--935, 2009.

\bibitem{prats0}
B.~Eynard and A.~Prats Ferrer.
\newblock 2-matrix versus complex matrix model, integrals over the unitary
  group as triangular integrals.
\newblock {\em Comm. Math. Phys.}, 264(1):115--144, 2006.

\bibitem{Samelson}
H.~Samelson.
\newblock {\em Notes on {L}ie algebras}.
\newblock Universitext. Springer-Verlag, New York, second edition, 1990.

\bibitem{Macdonald}
I.G. Macdonald.
\newblock The volume of a compact lie group.
\newblock {\em Inv. Math.}, 56(2):93--95, 1980.

\end{thebibliography}
\end{document}